# On the support of solutions to the generalized KdV equation

Carlos E. Kenig
Department of Mathematics
University of Chicago
Chicago, IL 60637, USA

Gustavo Ponce
Department of Mathematics
University of California
Santa Barbara, CA 93106, USA

and

Luis Vega
Departamento de Matematicas
Universidad del Pais Vasco
Apartado 644
48080 Bilbao, SPAIN.

1991 *Mathematics Subject Classification.* Primary 35Q53; secondary 35G25, 35D99.
C. E. K. and G. P. were supported by NSF grants. L. V. was supported by a DGICYT grant.

## §1. Introduction.

Consider the following question : Let $u = u(x,t)$ be a real valued solution of the $k$-generalized Korteweg-de Vries ($k$-gKdV) equation

$$(1.1) \qquad \partial_t u + \partial_x^3 u + u^k \partial_x u = 0, \qquad (x,t) \in \mathbb{R} \times (t_1, t_2), \quad k \in \mathbb{Z}^+,$$

with $t_1 < t_2$ which is sufficiently smooth and such that

$$(1.2) \qquad supp\ u(\cdot, t_j) \subseteq (a,b), \quad -\infty < a < b < \infty, \quad j=1,2.$$

Is $u \equiv 0$?

The first results in this direction are due to J.-C. Saut and B. Scheurer [SaSc]. They established the following unique continuation result.

**Theorem 1.1** [SaSc]. *Assume that $u = u(x,t)$ satisfies the equation*

$$(1.3) \qquad \partial_t u + \partial_x^3 u + \sum_{j=0}^{2} r(x,t)\, \partial_x^j u = 0, \qquad (x,t) \in (a,b) \times (t_1, t_2),$$

*with*

$$(1.4) \qquad r_j \in L^\infty((t_1, t_2) : L^2_{loc}((a,b))).$$

*If $u$ vanishes on an open set $\Omega \subseteq (a,b) \times (t_1, t_2)$, then $u$ vanishes in the horizontal components of $\Omega$, i.e. the set*

$$(1.5) \qquad \{(x,t) \in (a,b) \times (t_1, t_2) : \exists y\ s.t.\ (y,t) \in \Omega\}.$$

As a consequence they obtained the following result.

**Corollary 1.2** [SaSc]. *If $u$ is a sufficiently smooth solution of the equation (1.1) with*

$$(1.6) \qquad supp\ u(\cdot, t) \subseteq (a,b)^c, \quad \forall t \in (t_1, t_2),$$

*then* $u \equiv 0$.

The key step in Saut-Scheurer's argument is the following Carleman estimate :

Assume $(0,0) \in \Omega$ then $\exists \delta_0$, $M$, $K > 0$ such that

(1.7)
$$K \iint_\Omega |\partial_t u + \partial_x^3 u|^2 \exp(2\lambda\varphi) dx dt \geq \lambda \iint_\Omega |\partial_x^2 u|^2 \exp(2\lambda\varphi) dx dt$$
$$+ \lambda^2 \iint_\Omega |\partial_x u|^2 \exp(2\lambda\varphi) dx dt + \lambda^4 \iint_\Omega |u|^2 \exp(2\lambda\varphi) dx dt$$

for all $\lambda$ with $\lambda \delta \geq M$, $0 < \delta < \delta_0$ and $\varphi(x,t) = (x-\delta)^2 + \delta^2 t^2$.

In 1992, B. Zhang [Z] gave a positive answer to our question for the KdV equation

(1.8)
$$\partial_t u + \partial_x^3 u + u \partial_x u = 0$$

and for

(1.9)
$$\partial_t u + \partial_x^3 u - u^2 \partial_x u = 0,$$

using inverse scattering theory and Miura's transformation.

In 1997, J. Bourgain [B] used a different approach to reprove Corollary 1.2. His argument is based on the analyticity of the nonlinear term and the dispersion relation of the linear part of the equation. It also applies to higher order dispersive nonlinear models, and to higher spatial dimensions.

Recently, S. Tarama [T] showed that solutions $u(x,t)$ of the KdV equation (1.8) corresponding to data $u_0 \in L^2(\mathbb{R})$ such that

(1.10)
$$\int_{-\infty}^\infty (1+|x|)|u_0(x)| dx + \int_0^\infty e^{\delta |x|^{1/2}} |u_0(x)|^2 dx < \infty$$

for some $\delta > 0$, becomes analytic with respect to the space variable $x$ for $t > 0$. The proof is based on the inverse scattering method. Clearly this also provides a positive answer to our question in the case of the KdV equation.

The statement of our main result is the following.

**Theorem 1.3.** *Suppose that $u$ is a sufficiently smooth real valued solution of*

(1.11) $$\partial_t u + \partial_x^3 u + F(x, t, u, \partial_x u, \partial_x^2 u) = 0, \qquad (x, t) \in \mathbb{R} \times [t_1, t_2],$$

*where $F \in C_b^5$ in $(x,t)$, of polynomial growth in the other variables, at least quadratic in $u$, $\partial_x u$, $\partial_x^2 u$ in the terms involving $\partial_x^2 u$, i.e. $\partial F(x, t, 0, 0, 0)/\partial x_5 = 0$, for any $(x, t) \in \mathbb{R} \times [t_1, t_2]$.*

*Suppose also that $u$ has some decay if $F$ is just quadratic in $u$, $\partial_x u$, $\partial_x^2 u$ in the terms involving $\partial_x^2 u$.*

*If*

(1.12) $$\begin{aligned} \operatorname{supp} u(\cdot, t_j) &\subseteq (-\infty, b), \quad j = 1, 2, \\ (\operatorname{or}\ \operatorname{supp} u(\cdot, t_j) &\subseteq (a, \infty), \quad j = 1, 2), \end{aligned}$$

*then $u \equiv 0$.*

Remarks : (a) For the KdV equation (1.8) B. Zhang [Z] also had a similar result, i.e. one-sided support (1.12). Also for the KdV equation as a consequence of S. Tarama's result in [T] one finds that $u \equiv 0$, if there exists $t_1 < t_2$ such that $\operatorname{supp} u(\cdot, t_1) \subseteq (-\infty, b)$ and $\operatorname{supp} u(\cdot, t_2) \subseteq (a, \infty)$.

(b) It will be clear from our proof below that the result in Theorem 1.3 extends to complex valued solutions for the cases where energy estimates are available (see Lemma 2.1). For example, this holds for the equation

(1.13) $$\partial_t u + \partial_x^3 u \pm |u|^{2k} \partial_x u = 0, \quad k \in \mathbb{Z}^+.$$

(c) Although here we are not concerned with the minimal regularity assumptions on the solution $u$ required in Theorem 1.3, we remark that it suffices to assume that

(1.14) $$u \in C([t_1, t_2] : H^6(\mathbb{R}) \cap L^2(|x|^6 dx)) \cap C^1([t_1, t_2] : H^3(\mathbb{R})),$$

for the general case and

(1.15) $$u \in C([t_1, t_2] : H^6(\mathbb{R})) \cap C^1([t_1, t_2] : H^3(\mathbb{R})),$$

if the nonlinearity $F$ is at least cubic in $u$, $\partial_x u$, $\partial_x^2 u$ in the terms involving $\partial_x^2 u$, see [KePoVe3].

(d) To simplify the exposition we will carry out the details only in the case of the $k$-generalized KdV equation (1.1). In this case it suffices to assume that

$$(1.16) \qquad u \in C([t_1, t_2] : H^4(\mathbb{R})) \cap C^1([t_1, t_2] : H^1(\mathbb{R})).$$

For the existence theory we refer to [KePoVe2].

(e) Theorem 1.3 and its proof below extend to higher order dispersive models of the form

$$(1.17) \qquad \partial_t u + \partial_x^{2j+1} u + F(x, t, u, .., \partial_x^{2j} u) = 0, \quad j \in \mathbb{Z}^+,$$

whose local theory was developed in [KePoVe3].

(f) It should be remarked that we do not assume analyticity of the nonlinearity $F$.

The rest of the paper is organized as follows. In Section 2 we prove Theorem 1.3 assuming a key step in the proof, Lemma 2.3, whose proof is given in Section 3. Section 4 contains some remarks concerning the proofs and extensions of some of the results used in the proof of Theorem 1.3.

## §2. Proof of Theorem 1.3.

Without loss of generality we assume that $t_1 = 0$, $t_2 = 1$. Thus,

$$(2.1) \qquad supp \; u(\cdot, 0), \; supp \; u(\cdot, 1) \subseteq (-\infty, b).$$

We need some preliminary results.

The first one is concerned with the decay properties of solutions to the $k$-gKdV. The idea goes back to T. Kato [K].

**Lemma 2.1.** *Let $u = u(x, t)$ be a real valued solution of the $k$-gKdV equation (1.1) such that*

$$(2.2) \qquad \sup_{t \in [0,1]} \|u(\cdot, t)\|_{H^1} < \infty$$

and such that for a given $\beta > 0$

(2.3) $$e^{\beta x} u_0 \in L^2(\mathbb{R}).$$

Then

(2.4) $$e^{\beta x} u \in C([0,1] : L^2(\mathbb{R})).$$

*Proof.*

Let $\varphi_n \in C^\infty(\mathbb{R})$, with $\varphi_n(x) = e^{\beta x}$ for $x \leq n$, $\varphi_n(x) = e^{2\beta n}$ for $x > 10n$, $\varphi_n(x) \leq e^{\beta x}$, $0 \leq \varphi_n'(x) \leq \beta \varphi_n(x)$, and $|\varphi_n^{(j)}(x)| \leq \beta^j \varphi_n(x)$, $j = 2, 3$.

Multiplying the equation (1.1) by $u\varphi_n$, and integrating by parts we get

(2.5) $$\frac{1}{2}\frac{d}{dt}\int u^2 \varphi_n dx + \frac{3}{2}\int (\partial_x u)^2 \varphi_n' dx - \frac{1}{2}\int u^2 \varphi_n^{(3)} dx - \frac{1}{k+2}\int u^{k+2} \varphi_n' dx = 0.$$

Thus

(2.6) $$\frac{d}{dt}\int u^2 \varphi_n dx \leq \beta^3 \int u^2 \varphi_n dx + \frac{2\beta}{k+2}\|u\|_{L^\infty}^k \int u^2 \varphi_n dx$$

and

(2.7) $$\sup_{t\in[0,1]}\int u^2(x,t)\varphi_n(x)dx \leq \left(\int u_0^2 e^{\beta x} dx\right) exp(C^*),$$

where

(2.8) $$C^* = \beta^3 + \frac{2\beta}{k+2}\|u\|_{L^\infty(\mathbb{R}\times[0,1])}^k.$$

Now taking $n \uparrow \infty$ we obtain the desired result (2.4).

Lemma 2.1 has the following extension to higher derivatives.

**Lemma 2.2.** *Let $j \in \mathbb{Z}$, $j \geq 1$. Let $u = u(x,t)$ be a solution of the k-gKdV equation (1.1) such that*

$$\sup_{t \in [0,1]} \|u(\cdot, t)\|_{H^{j+1}} < \infty \tag{2.9}$$

*and for a given $\beta > 0$*

$$e^{\beta x} u_0, .., e^{\beta x} \partial_x^j u_0 \in L^2(\mathbb{R}). \tag{2.10}$$

*Then*

$$\sup_{t \in [0,1]} \|e^{\beta x} u(t)\|_{C^{j-1}} \leq c_j = c_j(u_0; C^*), \tag{2.11}$$

*with $C^*$ as in (2.8).*

Under the hypothesis (2.9)-(2.10) the result in [K] (Theorem 11.1) guarantees that $u \in C^\infty(\mathbb{R} \times (0,1])$.

To state the next results we need to introduce some notation,

$$f \in C^{3,1}(\mathbb{R}^2) \quad \text{if} \quad \partial_x f,\ \partial_x^2 f,\ \partial_x^3 f,\ \partial_t f \in C(\mathbb{R}^2), \tag{2.12}$$

and

$$f \in C_0^{3,1}(\mathbb{R}^2) \quad \text{if} \quad f \in C^{3,1}(\mathbb{R}^2) \text{ with compact support.} \tag{2.13}$$

Next, following the ideas in Kenig-Ruiz-Sogge [KeRuSo] and Kenig-Sogge [KeSo] we have the following Carleman Estimates.

**Lemma 2.3.** *If $f \in C_0^{3,1}(\mathbb{R}^2)$ (see (2.13)), then*

$$\|e^{\lambda x} f\|_{L^8(\mathbb{R}^2)} \leq c \|e^{\lambda x} \{\partial_t + \partial_x^3\} f\|_{L^{8/7}(\mathbb{R}^2)} \tag{2.14}$$

*for all $\lambda \in \mathbb{R}$, with $c$ independent of $\lambda$.*

The proof of Lemma 2.3, which is similar to those in [KeRuSo] and [KeSo], will be given in Section 3.

**Lemma 2.4.** *If $g \in C^{3,1}(\mathbb{R}^2)$ (see (2.12)) is such that*

(2.15) $$\text{supp } g \subseteq [-M, M] \times [0, 1]$$

*and*

(2.16) $$g(x, 0) = g(x, 1) = 0, \quad \forall x \in \mathbb{R},$$

*then*

(2.17) $$\|e^{\lambda x} g\|_{L^8(\mathbb{R}\times[0,1])} \leq c \|e^{\lambda x} \{\partial_t + \partial_x^3\} g\|_{L^{8/7}(\mathbb{R}\times[0,1])}$$

*for all $\lambda \in \mathbb{R}$, with $c$ independent of $\lambda$.*

*Proof.*

Let $\theta_\epsilon \in C_0^\infty(\mathbb{R})$, with $\theta_\epsilon(t) = 1$ for $t \in (\epsilon, 1-\epsilon)$, $0 \leq \theta_\epsilon(t) \leq 1$ and $|\theta_\epsilon'(t)| \leq c/\epsilon$.

Let

(2.18) $$g_\epsilon(x, t) = \theta_\epsilon(t) g(x, t).$$

we will apply (2.14) to $g_\epsilon$ for all $\epsilon > 0$. On the one hand

(2.19) $$\|e^{\lambda x} g_\epsilon\|_{L^8(\mathbb{R}^2)} = \|e^{\lambda x} g_\epsilon\|_{L^8(\mathbb{R}\times[0,1])} \to \|e^{\lambda x} g\|_{L^8(\mathbb{R}\times[0,1])} \text{ as } \epsilon \downarrow 0.$$

On the other hand,

(2.20) $$\{\partial_t + \partial_x^3\} g_\epsilon = \theta_\epsilon(t) \{\partial_t + \partial_x^3\} g + \theta_\epsilon'(t) g,$$

(2.21) $$\|e^{\lambda x} \theta_\epsilon(t) \{\partial_t + \partial_x^3\} g\|_{L^{8/7}(\mathbb{R}^2)} \to \|e^{\lambda x} \{\partial_t + \partial_x^3\} g\|_{L^{8/7}(\mathbb{R}\times[0,1])},$$

and

(2.22) $$\|\theta_\epsilon'(t) g\|_{L^{8/7}(\mathbb{R}^2)} \leq \frac{c}{\epsilon} \left( \int_0^\epsilon \int_{-M}^M |g(x,t)|^{8/7} dx dt \right)^{7/8}$$
$$+ \frac{c}{\epsilon} \left( \int_{1-\epsilon}^1 \int_{-M}^M |g(x,t)|^{8/7} dx dt \right)^{7/8} \to 0 \text{ as } \epsilon \downarrow 0.$$

from the mean value theorem.

**Lemma 2.5.** *Let $g \in C^{3,1}(\mathbb{R} \times [0,1])$ (see (2.12)). Suppose that*

$$\sum_{j \leq 2} |\partial_x^j g(x,t)| \leq c_\beta \, e^{-\beta |x|}, \quad t \in [0,1], \quad \forall \beta > 0, \tag{2.23}$$

*and*

$$g(x,0) = g(x,1) = 0, \quad \forall x \in \mathbb{R}. \tag{2.24}$$

*Then*

$$\|e^{\lambda x} g\|_{L^8(\mathbb{R} \times [0,1])} \leq c_0 \, \|e^{\lambda x} \{\partial_t + \partial_x^3\} g\|_{L^{8/7}(\mathbb{R} \times [0,1])} \tag{2.25}$$

*for all $\lambda \in \mathbb{R}$, with $c_0$ independent of $\lambda$.*

*Proof.*

Let $\phi \in C_0^\infty(\mathbb{R})$ be an even, nonincreasing function for $x > 0$ with $\phi(x) = 1$, $|x| \leq 1$, and $supp \, \phi \subseteq [-2, 2]$. Define $\phi_M(x) = \phi(x/M)$.

Let $g_M(x,t) = \phi_M(x) \, g(x,t)$.

Since

$$\begin{aligned} \{\partial_t + \partial_x^3\} g_M &= \phi_M \{\partial_t + \partial_x^3\} g + 3 \partial_x \phi_M \partial_x^2 g + 3 \partial_x^2 \phi_M \partial_x g + \partial_x^3 \phi_M g \\ &= \phi_M \{\partial_t + \partial_x^3\} g + E_1 + E_2 + E_3, \end{aligned} \tag{2.26}$$

applying Lemma 2.4 to $g_M(x,t)$ we get

$$\begin{aligned} \|e^{\lambda x} g_M\|_{L^8(\mathbb{R} \times [0,1])} &\leq c \, \|e^{\lambda x} \{\partial_t + \partial_x^3\} g_M\|_{L^{8/7}(\mathbb{R} \times [0,1])} \\ &\leq c \, \|e^{\lambda x} \phi_M \{\partial_t + \partial_x^3\} g\|_{L^{8/7}(\mathbb{R} \times [0,1])} + c \sum_{j=1}^3 \|e^{\lambda x} E_j\|_{L^{8/7}(\mathbb{R} \times [0,1])}. \end{aligned} \tag{2.27}$$

We need to show that the terms involving the $L^{8/7}$-norm of the "errors" $E_1$, $E_2$, and $E_3$ in (2.27) tend to zero as $M \uparrow \infty$. It suffices to consider one of them, say $E_1$, since the proof for $E_2$, $E_3$ is similar. Also it will be clear from the argument given below that it suffices to consider only the case $x > 0$ and $\lambda > 0$. From (2.23) with $\beta > \lambda$ it follows that

$$\|e^{\lambda x} E_1\|^{8/7}_{L^{8/7}(\mathbb{R}^+\times[0,1])} = 3^{8/7} \int_0^1 \int_M^{2M} |e^{\lambda x} \partial_x \phi_M \partial_x^2 g|^{8/7} dx dt$$

(2.28)
$$\leq c \int_0^1 \int_M^{2M} \left|\frac{e^{\lambda x}}{M} \partial_x^2 g\right|^{8/7} dx dt$$

$$\leq c \int_0^1 \int_M^{2M} e^{8\lambda x/7} e^{-8\beta x/7} dx dt \to 0 \quad \text{as} \quad M \uparrow \infty.$$

Thus, taking the limit as $M \uparrow \infty$ in (2.27) and using (2.28) we obtain (2.25).

**Lemma 2.6.** *Suppose $u = u(x,t) \in C([0,1] : H^4(\mathbb{R})) \cap C^1([0,1] : H^1(\mathbb{R}))$ satisfies the equation*

(2.29)
$$\partial_t u + \partial_x^3 u + u^k \partial_x u = 0, \quad (x,t) \in \mathbb{R} \times [0,1]$$

*with*

(2.30)
$$supp\, u(x,0) \subseteq (-\infty, b].$$

*Then for any $\beta > 0$*

(2.31)
$$\sum_{j\leq 2} |\partial_x^j u(x,t)| \leq c_{b,\beta}\, e^{-\beta x}, \quad \text{for } x > 0,\ t \in [0,1].$$

*Proof.*

It follows from Lemma 2.2.

*Proof of Theorem 1.3.*

We will show that there exists a large number $R > 0$ such that

(2.32)
$$supp\, u(\cdot,t) \subseteq (-\infty, 2R], \quad \forall t \in [0,1].$$

Then Saut-Schaurer's result (Theorem 1.1) completes the proof.

Let $\mu \in C^\infty(\mathbb{R})$ be a nondecreasing function such that $\mu(x) = 0$, $x \leq 1$ and $\mu(x) = 1$, $x \geq 2$. Let $\mu_R(x) = \mu(x/R)$.

Define

(2.33) $\quad V(x,t) = u^{k-1}(x,t)\partial_x u(x,t) \in L^p(\mathbb{R} \times [0,1])$, $\forall p \in [1,\infty]$ (by (2.31)),

and

(2.34) $$u_R(x,t) = \mu_R(x)u(x,t).$$

Combining our assumptions (see (1.16)) and Lemma 2.6 we can apply Lemma 2.5 to $u_R(x,t)$ for $R$ sufficiently large. Thus, using that

$$
\begin{aligned}
\{\partial_t + \partial_x^3\}u_R(x,t) &= \{\partial_t + \partial_x^3\}(\mu_R u) \\
&= \mu_R V u + 3\partial_x \mu_R \partial_x^2 u + 3\partial_x^2 \mu_R \partial_x u + \partial_x^3 \mu_R u \\
&= \mu_R V u + F_1 + F_2 + F_3 = \mu_R V u + F_R,
\end{aligned}
$$
(2.35)

it follows that

(2.36)
$$
\begin{aligned}
\|e^{\lambda x}\mu_R u\|_{L^8(\mathbb{R}\times[0,1])} &\leq c_0 \|e^{\lambda x}\{\partial_t + \partial_x^3\}(\mu_R u)\|_{L^{8/7}(\mathbb{R}\times[0,1])} \\
&\leq c_0 \|e^{\lambda x}\mu_R V u\|_{L^{8/7}(\mathbb{R}\times[0,1])} + c_0 \|e^{\lambda x} F_R\|_{L^{8/7}(\mathbb{R}\times[0,1])}.
\end{aligned}
$$

where $c_0$ is the constant coming from Lemma 2.5, (2.25). Then

(2.37)
$$
\begin{aligned}
c_0 \|e^{\lambda x}\mu_R V u\|_{L^{8/7}(\mathbb{R}\times[0,1])} \\
\leq c_0 \|e^{\lambda x}\mu_R u\|_{L^8(\mathbb{R}\times[0,1])} \|V\|_{L^{4/3}(\{x\geq R\}\times[0,1])}.
\end{aligned}
$$

Now we fix $R$ so large such that

(2.38) $$c_0 \|V\|_{L^{4/3}(\{x\geq R\}\times[0,1])} \leq 1/2.$$

From (2.36)-(2.38) one finds that

(2.39) $$\|e^{\lambda x}(\mu_R u)\|_{L^8(\mathbb{R}\times[0,1])} \leq 2c_0 \|e^{\lambda x} F_R\|_{L^{8/7}(\mathbb{R}\times[0,1])}.$$

As in the proof of Lemma 2.5 to estimate the left hand side of (2.39) it suffices to consider one of the terms in $F_R$, say $F_2$, since the proofs for $F_1$, $F_3$ are similar. We recall that the supports of the $F_j$'s are contained in the interval $[R, 2R]$. Thus,

$$
\begin{aligned}
(2.40) \quad 2c_0 \|e^{\lambda x} F_2\|_{L^{8/7}(\mathbb{R}\times[0,1])} &\leq \frac{2c_0}{R^2} \left( \int_0^1 \int_R^{2R} e^{8\lambda x/7} |\partial_x u(x,t)|^{8/7} dx dt \right)^{7/8} \\
&\leq \frac{2c_0}{R^2} e^{2\lambda R} \left( \int_0^1 \int_R^{2R} |\partial_x u(x,t)|^{8/7} dx dt \right)^{7/8}.
\end{aligned}
$$

On the other hand,

$$
(2.41) \quad \|e^{\lambda x}(\mu_R u)\|_{L^8(\mathbb{R}\times[0,1])} \geq \left( \int_0^1 \int_{x>2R} e^{8\lambda x} |u(x,t)|^8 dx dt \right)^{1/8}.
$$

Combining (2.39)-(2.41) we conclude that

$$
\begin{aligned}
(2.42) \quad &\left( \int_0^1 \int_{x>2R} e^{8\lambda(x-2R)} |u(x,t)|^8 dx dt \right)^{1/8} \\
&\leq \frac{2c_0}{R^2} \left( \int_0^1 \int_R^{2R} |\partial_x u(x,t)|^{8/7} dx dt \right)^{7/8}.
\end{aligned}
$$

Now letting $\lambda \uparrow \infty$ it follows that

$$(2.43) \quad u(x,t) \equiv 0 \quad \text{for} \quad x > 2R, \ t \in [0,1],$$

which yields the proof.

## §3. Proof of Lemma 2.3.

We shall prove that if $f \in C_0^{3,1}(\mathbb{R}^2)$, see (2.13), then

$$(3.1) \quad \|e^{\lambda x} f\|_{L^8(\mathbb{R}^2)} \leq c \|e^{\lambda x}\{\partial_t + \partial_x^3\} f\|_{L^{8/7}(\mathbb{R}^2)},$$

for all $\lambda \in \mathbb{R}$, with $c$ independent of $\lambda$.

We divide the proof into five steps.

STEP 1: It suffices to consider the cases $\lambda = \pm 1$ in (3.1).

*Proof.*

To prove the claim we observe that the case $\lambda = 0$ follows from the case $\lambda \neq 0$ by taking the limit as $\lambda \to 0$. So we can restrict ourselves to the case $\lambda \neq 0$.

Consider the case $\lambda > 0$ (the proof for $\lambda < 0$ is similar). Assume that

$$\|e^x f\|_{L^8(\mathbb{R}^2)} \leq c \|e^x \{\partial_t + \partial_x^3\} f\|_{L^{8/7}(\mathbb{R}^2)} \tag{3.2}$$

for all $f \in C_0^{3,1}(\mathbb{R}^2)$ with $c$ independent of $\lambda$.

Defining

$$f_\lambda(x, t) = f(x/\lambda, t/\lambda^3), \tag{3.3}$$

one has that

$$\{\partial_t + \partial_x^3\} f_\lambda(x, t) = \frac{1}{\lambda^3} \left( \partial_t f(x/\lambda, t/\lambda^3) + \partial_x^3 f(x/\lambda, t/\lambda^3) \right). \tag{3.4}$$

From the change of variables

$$(y, s) = (x/\lambda, t/\lambda^3), \qquad dx\, dt = \lambda^4\, dy\, ds, \tag{3.5}$$

it follows that

$$\|e^x f_\lambda\|_{L^8} = \lambda^{4/8} \|e^{\lambda y} f\|_{L^8} = \lambda^{1/2} \|e^{\lambda y} f\|_{L^8}, \tag{3.6}$$

and

$$\begin{aligned}\|e^x \{\partial_t + \partial_x^3\} f_\lambda\|_{L^{8/7}} &= \frac{\lambda^{4 \cdot 7/8}}{\lambda^3} \|e^{\lambda y} \{\partial_s + \partial_y^3\} f\|_{L^{8/7}} \\ &= \lambda^{1/2} \|e^{\lambda y} \{\partial_s + \partial_y^3\} f\|_{L^{8/7}}\end{aligned} \tag{3.7}$$

Inserting (3.6)-(3.7) into (3.2) we obtain (3.1), which proves the claim.

STEP 2: To prove (3.2) it suffices to establish the following inequality

$$\|g\|_{L^8} \leq c \|\{\partial_t + \partial_x^3 - 3\partial_x^2 + 3\partial_x - 1\} g\|_{L^{8/7}}, \tag{3.8}$$

for any $g \in C_0^{3,1}(\mathbb{R}^2)$, see (2.13).

*Proof.*

Let

$$(3.9) \qquad g(x,t) = e^x f(x,t).$$

Since

$$(3.10) \qquad e^x \{\partial_t + \partial_x^3\} f = \{\partial_t + \partial_x^3 - 3\partial_x^2 + 3\partial_x - 1\} g,$$

we obtain (3.2).

<u>STEP 3</u>: It suffices to prove the inequality (3.8) without the term in the left hand side involving the derivatives of order 1 in the $x$-variable. In other words, to prove (3.8) it suffices to show

$$(3.11) \qquad \|h\|_{L^8} \le c\|\{\partial_t + \partial_x^3 - 3\partial_x^2 - 1\} h\|_{L^{8/7}},$$

for any $h \in C_0^{3,1}(\mathbb{R}^2)$, see (2.13).

*Proof.*

Using the change of variables

$$(3.12) \qquad y = x/3 + t, \quad s = t, \quad (x = 3(s-y), \quad t = s), \quad (dy\,ds = dx\,dt/3),$$

and the notation

$$(3.13) \qquad h(y,s) = g(x,t)$$

it follows that

$$(3.14) \qquad \frac{\partial h}{\partial y} = -3 \frac{\partial g}{\partial x}, \qquad \frac{\partial h}{\partial s} = 3 \frac{\partial g}{\partial x} + \frac{\partial g}{\partial t}.$$

Thus, (3.8) can be written in the equivalent form

$$(3.15) \qquad \|h\|_{L^8} \le c \|\{\partial_s - \frac{1}{27} \partial_y^3 - \frac{1}{3} \partial_y^2 - 1\} h\|_{L^{8/7}}.$$

Finally, making another change of variables

$$(3.16) \quad z = -\frac{1}{3}y, \quad t = s, \quad (\partial_z^3 = -\frac{1}{27}\partial_y^3, \ \partial_z^2 = \frac{1}{9}\partial_y^2, \ \partial_t = \partial_s),$$

it follows that (3.15) is equivalent to

$$(3.17) \quad \|h\|_{L^8} \leq c\|\{\partial_t + \partial_z^3 - 3\partial_z^2 - 1\}h\|_{L^{8/7}},$$

which proves the claim.

<u>STEP 4</u> We will need the following results (Lemmas 3.1-3.2). The first one is an estimate of Strichartz type.

**Lemma 3.1.**

$$(3.18) \quad \|\int_{\mathbb{R}} e^{i(x,t)\cdot(\xi,\xi^3)} \hat{f}(\xi,\xi^3)\,d\xi\|_{L^8(\mathbb{R}^2)} \leq c\|f\|_{L^{8/7}(\mathbb{R}^2)}.$$

where $\hat{}$ denotes the Fourier transform.

*Proof of Lemma 3.1.*

Using the notation

$$(3.19) \quad U(t)v_0(x) = \int_{-\infty}^{\infty} e^{i(t\xi^3 + x\xi)} \hat{v}_0(\xi)d\xi = (e^{it\xi^3}\hat{v}_0)^{\vee}(x,t),$$

the inequality (3.18) can be written as

$$(3.20) \quad \|\int_{-\infty}^{\infty} U(t-t')f(\cdot,t')dt'\|_{L^8(\mathbb{R}^2)} \leq c\|f\|_{L^{8/7}(\mathbb{R}^2)}$$

whose proof can be found in [GiTs], (Lemma 2.1) or in [KePoVe1], (Theorem 2.1).

**Lemma 3.2.**

$$(3.21) \quad \|h\|_{L^8} \leq c\|\{\partial_t + \partial_x^3 + a\}h\|_{L^{8/7}},$$

for any $h \in C_0^{3,1}(\mathbb{R}^2)$, see (2.13), with c independent of $a \in \mathbb{R}$.

*Proof of Lemma 3.2.*

Using the notation introduced in (3.19) we recall the decay estimate

$$\text{(3.22)} \qquad \|U(t)v_0\|_{L^8(\mathbb{R})} \leq \frac{c}{|t|^{1/4}} \|v_0\|_{L^{8/7}(\mathbb{R})},$$

which follows by interpolating the estimates

$$\text{(3.23)} \qquad \|U(t)v_0\|_{L^2} = \|v_0\|_{L^2}, \qquad \|U(t)v_0\|_{L^\infty} \leq \frac{c}{|t|^{1/3}} \|v_0\|_{L^1}.$$

An homogeneity argument, similar to that given in Step 1, shows that it suffices to consider only the case $|a| = 1$. We thus need to prove the multiplier estimate

$$\text{(3.24)} \qquad \left\| \left( \frac{1}{i(\tau - \xi^3) \mp 1} \hat{h}(\xi, \tau) \right)^{\vee} \right\|_{L^8} = \left\| \left( \frac{1}{\tau - \xi^3 \pm i} \hat{h}(\xi, \tau) \right)^{\vee} \right\|_{L^8} \leq c \|h\|_{L^{8/7}}.$$

Let $S_\pm$ denotes the operator

$$\text{(3.25)} \qquad S_\pm h(x,t) = \int_{-\infty}^{\infty} \int_{-\infty}^{\infty} e^{i(x,t)\cdot(\xi,\tau)} \frac{1}{\tau - \xi^3 \pm i} \hat{h}(\xi, \tau) \, d\tau \, d\xi.$$

Let

$$\text{(3.26)} \qquad b_\pm(s) = \int_{-\infty}^{\infty} e^{i\tau s} \frac{1}{\tau \pm i} \, d\tau$$

so that

$$\text{(3.27)} \qquad S_\pm h(x,t) = \int_{-\infty}^{\infty} \left( \int_{-\infty}^{\infty} (h(\cdot, t-s))^{\wedge}(\xi) \, e^{is\xi^3} \, e^{ix\xi} \, d\xi \right) b_\pm(s) ds.$$

Thus,

$$\text{(3.28)} \qquad S_\pm h(x,t) = \int_{-\infty}^{\infty} U(s) h(\cdot, t-s) \, b_\pm(s) ds.$$

Note that

$$\text{(3.29)} \qquad \|b_\pm\|_{L^\infty} \leq c,$$

which combined with (3.22) leads to

$$\text{(3.30)} \qquad \|S_\pm h(\cdot, t)\|_{L^8(\mathbb{R})} \leq c \int_{-\infty}^{\infty} \|h(\cdot, t-s)\|_{L^{8/7}} \frac{ds}{|s|^{1/4}}.$$

Now $1/8 = 7/8 - 3/4$, and so fractional integration completes the proof.

STEP 5 To complete the proof of Lemma 2.3 we just need to prove (3.11), i.e.

$$\|h\|_{L^8} \leq c \|\{\partial_t + \partial_x^3 - 3\partial_x^2 - 1\}h\|_{L^{8/7}} \tag{3.31}$$

for any $h \in C_0^{3,1}(\mathbb{R}^2)$, see (2.13).

Taking Fourier transform, in space and time variables, in the left hand side of (3.31) we get

$$[i\tau - i\xi^3 + 3\xi^2 - 1]\hat{h}(\xi, \tau). \tag{3.32}$$

We consider the pair of points

$$P_\pm = (\xi_0^\pm, \tau_0^\pm) = \pm \left(\frac{1}{\sqrt{3}}, \left(\frac{1}{\sqrt{3}}\right)^3\right) \tag{3.33}$$

where the symbol in (3.32) vanishes. We recall that $h$ has compact support so its Fourier transform has an analytic continuation to $\mathbb{C}^2$. Hence, it suffices to prove (3.31) for any $h \in \mathcal{S}(\mathbb{R}^2)$ with $\hat{h}$ vanishing at $P_\pm$.

So we are then reduced to showing the multiplier inequality

$$\|\mathcal{M}h\|_{L^8(\mathbb{R}^2)} = \left\| \left( \frac{1}{i(\tau - \xi^3) + 3\xi^2 - 1} \hat{h} \right)^\vee \right\|_{L^8(\mathbb{R}^2)} \leq c\|h\|_{L^{8/7}(\mathbb{R}^2)}, \tag{3.34}$$

for such $h$'s.

It suffices to prove (3.34) assuming that

$$\operatorname{supp} \hat{h} \subseteq \{(\xi, \tau) : \xi \geq 0\}, \tag{3.35}$$

since the proof for the case

$$\operatorname{supp} \hat{h} \subseteq \{(\xi, \tau) : \xi < 0\}, \tag{3.36}$$

is similar.

We now recall a variant of Littlewood-Paley theory. Let

$$\widehat{L_k f}(\xi, \tau) = \chi_{[1/2, 1]}(|\xi - \xi_0^+|/2^{-k}) \hat{f}(\xi, \tau), \tag{3.37}$$

where $k \in \mathbb{Z}$ and $\chi_A(\cdot)$ is the characteristic function of the set $A$. Then for each $p \in (1, \infty)$ we have

$$\|f\|_{L^p(\mathbb{R}^2)} \simeq \left\| \left( \sum_{k \in \mathbb{Z}} |L_k f|^2 \right)^{1/2} \right\|_{L^p(\mathbb{R}^2)}. \tag{3.38}$$

Thus it suffices to establish (3.34) for each $L_k h$ with a constant independent of $k$, since using Minkowski's integral inequality ($8/7 < 2 < 8$) one has that

$$\|\mathcal{M} h\|_{L^8} \simeq \left\| \left( \sum_{k \in \mathbb{Z}^+} |L_k(\mathcal{M} h)|^2 \right)^{1/2} \right\|_{L^8}$$

$$= \left\| \left( \sum_{k \in \mathbb{Z}^+} |\mathcal{M}(L_k h)|^2 \right)^{1/2} \right\|_{L^8} \leq \left( \sum_{k \in \mathbb{Z}^+} \|\mathcal{M}(L_k h)\|_{L^8}^2 \right)^{1/2} \tag{3.39}$$

$$\leq c \left( \sum_{k \in \mathbb{Z}^+} \|L_k h\|_{L^{8/7}}^2 \right)^{1/2} \leq c \left\| \left( \sum_{k \in \mathbb{Z}^+} |L_k h|^2 \right)^{1/2} \right\|_{L^{8/7}}$$

$$\leq c \|h\|_{L^{/8/7}}.$$

Therefore, we shall prove the multiplier estimate (3.34) when

$$supp \, \hat{h} \subseteq \{(\xi, \tau) : \xi \geq 0, \ 2^{-k-1} \leq |\xi - \xi_0^+| \leq 2^{-k}\}, \tag{3.40}$$

We split the proof of (3.40) in two cases.

<u>CASE 1</u> : $k \leq 0$.

In this case, if $\xi \in supp \, \hat{h}$ then

$$|3\xi^2 - 1| \simeq |\xi - \xi_0^+| |\xi + \xi_0^+| \simeq 2^{-k} \cdot 2^{-k}. \tag{3.41}$$

Using Lemma 3.1 we just need to bound the multiplier

$$\frac{1}{i(\tau - \xi^3) + 3\xi^2 - 1} - \frac{1}{i(\tau - \xi^3) + 2^{-2k}}$$

(3.42)
$$= \frac{2^{-2k} - (3\xi^2 - 1)}{(i(\tau - \xi^3) + 3\xi^2 - 1)(i(\tau - \xi^3) + 2^{-2k})}.$$

Using the change of variables $\tau = \lambda + \xi^3$ write

$$\int_{-\infty}^{\infty}\int_{-\infty}^{\infty} \frac{e^{i(x,t)\cdot(\xi,\tau)}\left(2^{-2k} - (3\xi^2 - 1)\right)}{(i(\tau - \xi^3) + 3\xi^2 - 1)(i(\tau - \xi^3) + 2^{-2k})} \hat{h}(\xi,\tau) d\xi d\tau$$

$$= \int_{-\infty}^{\infty}\int_{-\infty}^{\infty} \frac{e^{i(x,t)\cdot(\xi,\lambda+\xi^3)}\left(2^{-2k} - (3\xi^2 - 1)\right)}{(i\lambda + 3\xi^2 - 1)(i\lambda + 2^{-2k})} \hat{h}(\xi, \lambda + \xi^3) d\xi d\lambda$$

(3.43)
$$= \int_{-\infty}^{\infty} e^{i\lambda t} \left( \int_{-\infty}^{\infty} \frac{e^{i(x,t)\cdot(\xi,\xi^3)}\left(2^{-2k} - (3\xi^2 - 1)\right)}{(i\lambda + 3\xi^2 - 1)(i\lambda + 2^{-2k})} \hat{h}(\xi, \lambda + \xi^3) d\xi \right) d\lambda$$

$$\equiv \Phi(x,t).$$

Defining

(3.44) $$\hat{h}_\lambda(\xi, \tau) = \frac{2^{-2k} - (3\xi^2 - 1)}{(i\lambda + 3\xi^2 - 1)(i\lambda + 2^{-2k})} \hat{h}(\xi, \tau + \lambda).$$

and using Lemma 3.1 and Minkowski's integral inequality we get

(3.45) $$\|\Phi\|_{L^8(\mathbb{R}^2)} \leq \int_{-\infty}^{\infty} \|h_\lambda\|_{L^{8/7}(\mathbb{R}^2)} d\lambda.$$

Now for $\lambda$ and $k$ fixed we consider the multiplier in (3.44) in the variable $\xi$, with $\xi > 0$, $|\xi - \xi_0^+| \simeq 2^{-k}$, $|3\xi^2 - 1| \simeq 2^{-2k}$, which has norm bounded by

(3.46) $$c \frac{2^{-2k}}{|\lambda|^2 + 2^{-4k}}.$$

Hence,

(3.47) $$\int_{-\infty}^{\infty} \|h_\lambda\|_{L^{8/7}(\mathbb{R}^2)} d\lambda \leq c \int_{-\infty}^{\infty} \frac{2^{-2k}}{|\lambda|^2 + 2^{-4k}} \|\tilde{h}_\lambda\|_{L^{8/7}(\mathbb{R}^2)} d\lambda \leq \|h\|_{L^{8/7}(\mathbb{R}^2)},$$

since $\tilde{h}_\lambda(x,t) = e^{-it\lambda}h(x,t)$, which combined with (3.45) yields the proof of case 1, i.e. $k \leq 0$.

CASE 2 : $k > 0$.

In this case, if $\xi \in \text{supp}\,\hat{h}$ then

(3.48) $$|3\xi^2 - 1| \simeq |\xi - \xi_0^+|\,|\xi + \xi_0^+| \simeq |\xi - \xi_0^+| \simeq 2^{-k}.$$

In this case we use Lemma 3.2 to substract

(3.49) $$\frac{1}{i(\tau - \xi^3) + 2^{-k}}$$

and argue exactly as before. The corresponding multiplier to (3.44) in this case is

(3.50) $$\frac{2^{-k} + 1 - 3\xi^2}{(i\lambda + 3\xi^2 - 1)(i\lambda + 2^{-k})}$$

which has norm bounded by

(3.51) $$c\,\frac{2^{-k}}{|\lambda|^2 + 2^{-2k}}.$$

## §4. Further Results.

In this section we extend Lemma 2.3 to higher order operators of the form considered in (1.17).

**Lemma 4.1.** *If $f \in C_0^\infty(\mathbb{R}^2)$, then for $j \in \mathbb{Z}^+$*

(4.1) $$\|e^{\lambda x} f\|_{L_t^q L_x^p} \leq c\|e^{\lambda x}\{\partial_t + \partial_x^{2j+1}\}f\|_{L_t^{q'} L_x^{p'}},$$

*for all $\lambda \in \mathbb{R}$, with $c$ independent of $\lambda \in \mathbb{R}$,*

(4.2) $$\frac{1}{2j+1} = \frac{2}{q} + \frac{2}{(2j+1)p}, \quad p \geq 2,$$

*and where $(p, p')$, $(q, q')$ are dual exponents, i.e. $1/p + 1/p' = 1$.*

*Proof.*

From the steps 1 and 2 in the proof of Lemma 2.3 one has that to establish (4.2) it suffices to show that

$$\|f\|_{L_t^q L_x^p} \leq c \, \|\{\partial_t + P(\partial_x)\}f\|_{L_t^{q'} L_x^{p'}}, \tag{4.3}$$

where

$$P(z) = (z-1)^{2j+1}. \tag{4.4}$$

Define $g = g(x,t)$ as

$$\{\partial_t + P(\partial_x)\}v(x,t) = g(x,t). \tag{4.5}$$

Taking Fourier transform in the $x$-variable in (4.5) we get

$$\partial_t \hat{v}(\xi,t) + P(i\xi)\,\hat{v}(\xi,t) = \hat{g}(\xi,t). \tag{4.6}$$

Since $v$ has compact support we conclude that $\hat{v}(\xi,T) = 0$ for any $T$ with $|T|$ large enough. Thus, from (4.6) it follows that

$$\partial_t(e^{t\,P(i\xi)}\,\hat{v})(\xi,t) = e^{t\,P(i\xi)}\,\hat{g}(\xi,t). \tag{4.7}$$

and

$$\begin{aligned}\hat{v}(\xi,t) &= \chi_{\{\xi:\operatorname{Re}(P(i\xi))\geq 0\}}(\xi) \int_{-\infty}^{t} e^{-(t-s)P(i\xi)}\,\hat{g}(\xi,s)ds \\ &\quad - \chi_{\{\xi:\operatorname{Re}(P(i\xi))<0\}}(\xi) \int_{t}^{\infty} e^{-(t-s)P(i\xi)}\,\hat{g}(\xi,s)ds.\end{aligned} \tag{4.8}$$

Hence,

$$\begin{aligned}v(x,t) &= \int_{-\infty}^{t}\int_{-\infty}^{\infty} K_+(x-y,t-s)g(y,s)dyds \\ &\quad - \int_{t}^{\infty}\int_{-\infty}^{\infty} K_-(x-y,t-s)g(y,s)dyds \\ &= L_+ g - L_- g,\end{aligned} \tag{4.9}$$

with

$$K_+(x,t) = \chi_{\{t:t\geq 0\}}(t) \int_{-\infty}^{\infty} e^{ix\cdot\xi - tP(i\xi)}\, \chi_{\{\xi:\operatorname{Re}(P(i\xi))\geq 0\}}(\xi)\, d\xi \tag{4.10}$$

and

$$K_-(x,t) = \chi_{\{t:t<0\}}(t) \int_{-\infty}^{\infty} e^{ix\cdot\xi - tP(i\xi)}\, \chi_{\{\xi:\operatorname{Re}(P(i\xi))<0\}}(\xi)\, d\xi. \tag{4.11}$$

We need the following results.

**Lemma 4.2.** *There exists a constant $c = c(j)$ such that*

(4.12) $$|K_\pm(x,t)| \leq \frac{c}{|t|^{1/(2j+1)}}.$$

*Proof.*

We consider the case of $K_+$. In this case the oscillatory part of the integral in (4.10) is given by the phase function $\phi(\xi) = Im\,(P(i\xi))$. We observe that

(4.13) $$\phi^{(2j+1)}(\xi) = (-1)^j\,(2j+1)!,$$

and

(4.14) $$\int_{\{\xi:Re\,(P(i\xi))\geq 0\}} |\partial_\xi(e^{-t\,P(i\xi)})|d\xi < c_j,$$

where $c_j$ depends on the numbers of changes of sign of $Re\,(P(i\xi))$. Hence, the proof of (4.12) follows from Van der Corput's lemma (see [S], Corollary in page 334).

**Lemma 4.3.** *For each $p \geq 2$ there exists a constant $c$ such that*

(4.15) $$\|\int_{-\infty}^{\infty} K_\pm(x-y,t)g(y)dy\|_{L_x^p} \leq \frac{c}{|t|^{(1/(2j+1))(1/p-1/p')}}\|g\|_{L_x^{p'}}$$

*Proof.*

For $p = 2$ we use Plancherel Theorem to get

(4.16) $$\|\int_{-\infty}^{\infty} K_+(x-y,t)g(y)dy\|_{L_x^2}^2$$
$$= \int_{-\infty}^{\infty} \chi_{\{\xi:Re\,(P(i\xi))\geq 0\}}(\xi)\,|e^{-t\,ReP(i\xi)}\,\hat{g}(\xi)|^2 d\xi \leq \|g\|_{L^2}.$$

The case $p = \infty$ follows from Lemma 4.2. Using the Riesz-Thorin theorem one extends the result to $p \in (2, \infty)$.

Finally the proof of Lemma 4.1 follows by combining (4.15), Minkowski's integral inequality, and Hardy-Littlewood-Sobolev inequality.